\documentclass[12pt]{amsart}
\usepackage{fullpage,url}
\usepackage{amssymb}
\usepackage[all]{xy} 

\DeclareFontEncoding{OT2}{}{} % to enable usage of cyrillic fonts

\newtheorem{theorem}{Theorem}[section]
\newtheorem{lemma}[theorem]{Lemma}

\newtheorem{proposition}[theorem]{Proposition}

\theoremstyle{definition}

\theoremstyle{remark}

\begin{document}

\title[paper]{Del Pezzo surfaces of degree 6}
%\subjclass[2000]{Primary 9999; Secondary 9999}
%\keywords{9999}
\author{Patrick Corn}
\address{Department of Mathematics, University of California, 
	Berkeley, CA 94720-3840, USA}
\email{corn@math.berkeley.edu}
\urladdr{http://math.berkeley.edu/\~{}corn}
\date{August 3, 2004}

\begin{abstract}

We give a correspondence which associates, to each Del Pezzo surface $X$ of degree $6$ 
over a field $k$ of characteristic $0$, a collection of data consisting of a 
Severi-Brauer variety/$k$ and a set of points defined over some extension of $k$.

\end{abstract}

\maketitle

The main results in this paper, and specifically Theorem \ref{theorem4}, give a way to 
describe Del Pezzo surfaces of degree $6$ over a field $k$ of characteristic $0$, via a 
correspondence with objects (Severi-Brauer varieties) which can be understood in a 
completely explicit way if $k$ is sufficiently nice (e.g. $k$ a number field).  

\section{Preliminaries}

In this paper, we will deal with varieties $V$ over a field $k$ of characteristic
$0$. If $L/k$ is
a field extension, then we write $V_L$ for the base extension $V \times_{{\rm Spec} 
\, k} {\rm Spec} \, L$, and $\overline V$ for $V_{\overline k}$.

A Del Pezzo surface over a number field $k$ is a smooth rational surface $X$ whose 
anticanonical sheaf $\omega_X^{-1}$ is ample. To each Del Pezzo surface $X$ is 
associated a number $d = (\omega_X,\omega_X)$ (where $(,)$ denotes intersection 
number), called the degree of $X$.

The results we need about Del Pezzo surfaces are summarized in the following 
proposition. We refer the interested reader to \cite{manin1974} for proofs and more 
details.

\begin{proposition} \label{prop1} 
Let $V$ be a Del Pezzo surface of degree $d$ over a field $k$.

(a) $1 \le d \le 9$.

(b) Pic $\overline V$ is a free abelian group of rank $10-d$.

(c) If $V' \to V$ is a birational morphism and $V'$ is a Del Pezzo surface, then 
$V$ is a Del Pezzo surface.

(d) Either $\overline V$ is isomorphic to the blowup of 
${\mathbb P}^2_{\overline k}$ at $r=9-d$ points $\{ x_1, \ldots, x_r \}$ in general
position, or $d=8$ and ${\overline V} \cong {\mathbb P}^1_{\overline k} \times 
{\mathbb P}^1_{\overline k}$. Conversely, if $r \le 6$, then any surface satisfying
this condition is a Del Pezzo surface of degree $d=9-r$. (For a set of $r \le 6$ 
points, ``general position" means that no three are collinear and no six lie on a 
conic.)

(e) Let $C$ be an exceptional curve; that is, $C$ is a curve on $\overline V$ 
such that $(C,C)=-1$ and $C \cong {\mathbb P}^1_{\overline k}$. Then if $r \le 6$,
the image of $C$ in ${\mathbb P}^2_{\overline k}$ under the isomorphism of (d) is
either: one of the $x_i$, a line passing through two of the $x_i$, or a conic
passing through five of the $x_i$. Conversely, each point, line, and conic in this
list gives rise to exactly one exceptional curve $C$.

\end{proposition}

{\em Proof:} These are (respectively) Theorem 24.3(i), Lemma 24.3.1, Corollary 
24.5.2, Theorem 24, and Theorem 26.2 of \cite{manin1974}. $\blacksquare$

\smallskip

The assumption that $r \le 6$ was made only to simplify the statements of (d) and 
(e); we will not be concerned with Del Pezzo surfaces of degree $1$ or $2$ in this 
paper.

\section{Severi-Brauer varieties: the basic construction}

If $V$ is a Del Pezzo surface over a field $k$, it is clear from the definition 
above that the exceptional curves on $\overline V$ are preserved by the action 
of $G_k := {\rm Gal}({\overline k}/k)$; this information can be very useful in 
investigating properties of these surfaces.

Now let $D$ be a Del Pezzo surface of degree 6. Let $E_1, E_2$, and $E_3$ be the 
exceptional curves corresponding to the blow-ups of the three points $x_1,x_2,x_3
\in {\mathbb P}^2_{\overline k}$ as in Proposition \ref{prop1}(d). Let $F_{12}$ be the 
exceptional curve corresponding to the line between $x_1$ and $x_2$, and define 
$F_{13}$ and $F_{23}$ similarly. By Proposition \ref{prop1}(e), the set $\{ E_1, E_2, 
E_3, F_{12}, F_{13}, F_{23} \}$ is precisely the set of exceptional curves on $D$.

We can now examine the possibilities for blowing down these curves to obtain other,
possibly simpler surfaces.

\begin{proposition} \label{prop2} Let $D$ be a Del Pezzo surface of degree six over a 
field $k$ of characteristic $0$. There is a field $L$ such that $[L:k] = 1$ or $2$ and 
surfaces $X$ and $Y$ defined over $L$ such that the triple $(D_L,X,Y)$ satisfies the 
following conditions:

(i) there is a morphism $\pi_X : D_L \to X$ which exhibits $D_L$ as the blow-up of 
$X$ at a $G_L$-stable set of three non-collinear points $\{P_1,P_2,P_3\} \in 
X({\overline k})$

(ii) there is a morphism $\pi_Y : D_L \to Y$ which exhibits $D_L$ as the blow-up
of $Y$ at a $G_L$-stable set of three non-collinear points $\{Q_1,Q_2,Q_3\} \in 
Y({\overline k})$

(iii) $\{ \pi_X^{-1}(P_i) : 1 \le i \le 3 \} \cup \{ \pi_Y^{-1}(Q_i) : 1 \le i 
\le 3 \}$ is a full set of six exceptional curves on $\overline D$.

(iv) $X$ and $Y$ are Severi-Brauer varieties of dimension $2$.

\end{proposition}

{\em Proof:} Let $L$ be the minimal field such that the sets $\{ E_1, E_2, E_3 \}$
and $\{F_{12}, F_{13}, F_{23} \}$ are both $G_L$-stable. Any element of $G_k$ 
either fixes both sets or switches them, so $L$ is either equal to $k$ or quadratic
over $k$. Let $X'$ and $Y'$ be the varieties obtained from blowing down 
$\{ E_1, E_2, E_3 \}$ and $\{ F_{12}, F_{13}, F_{23} \}$, respectively, over 
$\overline k$. Then they can naturally be descended to varieties $X$ and $Y$ 
defined over $L$ (see \cite{weil1956} for details on descent). Properties (i)-(iii) are 
immediate. 

To see property (iv), note that $X$ and $Y$ are Del Pezzo surfaces by Proposition 
\ref{prop1}(c). Now note that rank Pic $D = 4$ by Proposition \ref{prop1}(b), and 
blowing up at a point 
increases the rank of the Picard group by $1$, so rank Pic $\overline X$ must be $1$.
Then by Proposition \ref{prop1}(b) the degree of $X$ is $9$, which means $r=0$, so $X$ 
is a twist of ${\mathbb P}^2$. The same holds for $Y$. $\blacksquare$

\smallskip

Now we can also turn Proposition \ref{prop2} around: 

\begin{proposition} \label{prop3} Let $X$ be a Severi-Brauer variety over a field $L$ 
and let $\{ P_1,P_2,P_3 \}$ be a $G_L$-stable set of non-collinear points in 
$X({\overline L})$. Then there exist $S$ and $Y$ defined over $L$ such that the 
triple $(S,X,Y)$ satisfies conditions (i)-(iv) of Proposition 2.

\end{proposition}

{\em Proof:} To obtain $S$, simply blow up $X$ over $L$ at the given set of points.
To obtain $Y$, note that the three exceptional curves on $\overline S$ which 
are the inverse images of $\{ P_1,P_2,P_3 \}$ form a $G_L$-stable set (call it 
$C_1$), and since the full set $C$ of exceptional curves is $G_L$-stable, 
the complement $C \setminus C_1$ is also $G_L$-stable and can be blown down over
$L$ to obtain $Y$. Conditions (i)-(iii) are all obvious from the construction.
$\blacksquare$

\smallskip

There is a natural one-to-one correspondence between (isomorphism classes of) 
Severi-Brauer varieties over $L$ of dimension $2$ and (isomorphism classes of) central 
simple algebras over $L$ of dimension $9$; both are parametrized by 
$H^1(G_L,PGL_3({\overline L}))$. (Cf. \cite{serrelocalfields}.)

Given the result of proposition \ref{prop3}, the natural question to ask is: how are 
the central simple algebras corresponding to $X$ and $Y$ related? The answer is our
first main result. 

\begin{theorem} \label{theorem1} Let $X$ be a Severi-Brauer variety of dimension $2$ 
over a field $L$, equipped with a $G_L$-stable set $\{ P_1, P_2, P_3 \}$ of three 
non-collinear points. Construct $S$ and $Y$ as in Proposition \ref{prop2}. Let $x$ 
and $y$ be the central simple algebras corresponding to $X$ and $Y$ 
respectively. Then $y=x^{\rm op}$, the opposite algebra of $x$.

\end{theorem}

{\em Proof:} We simply unravel the definition of the correspondence between 
Severi-Brauer varieties and central simple algebras. First, choose an ordering
of the points $\{ P_1, P_2, P_3 \}$ and the points $\{ Q_1, Q_2, Q_3 \}$ of 
proposition 2 so that $\pi_X^{-1}(P_i) \cap \pi_Y^{-1}(Q_i) = \emptyset$ for all 
$i$.

Let $M$ be the minimal Galois extension of $L$ over which the $P_i$ are each 
individually defined. Then the $Q_i$ are all defined over $M$ as well. Also, 
$S_M$ is isomorphic to the blowup of ${\mathbb P}^2_M$ at the $P_i$, so we can choose
a point $P \in S(M)$ which lies over a point in ${\mathbb P}^2(M) \setminus \{ P_1, P_2,
P_3 \}$, so that $P$ does not lie on any exceptional curve. Now 
$X_M \cong {\mathbb P}^2_M$, and since the automorphism group of ${\mathbb 
P}^2_M$ acts transitively on sets of four $M$-points in general position, we can 
construct an isomorphism $\phi : X_M \to {\mathbb P}^2_M$ sending the 
points $P_1, P_2, P_3, (\pi_X)_M(P)$ on $X$ 
(notation as in Proposition \ref{prop2}) to $(1:0:0)$, $(0:1:0)$, $(0:0:1)$, and 
$(1:1:1)$, respectively. We can also construct an isomorphism 
$\psi : Y_M \to {\mathbb P}^2_M$ sending $Q_1, Q_2, Q_3, (\pi_Y)_M(P)$ on $Y$
(notation as in Proposition \ref{prop2}) to $(1:0:0), (0:1:0), (0:0:1)$, and $(1:1:1)$, 
respectively.

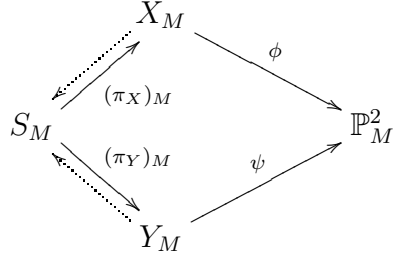
\[ \xymatrix{ & X_M \ar@<-0.5ex>@{.>}[dl] \ar[drr]^{\phi} & & \\
S_M \ar@<-0.5ex>[ur]_{(\pi_X)_M} \ar@<0.5ex>[dr]^{(\pi_Y)_M} & & & {\mathbb P}^2_M
\\
& Y_M \ar@<0.5ex>@{.>}[ul] \ar[urr]^{\psi} & &
} \]

Starting at ${\mathbb P}^2_M$ and going around counterclockwise in the diagram, we 
have a rational map $b := \psi \circ (\pi_Y)_M \circ (\pi_X)_M^{-1} \circ 
\phi^{-1}$. We can easily write down a formula for this map, as in 
\cite{hartshorne1977}, pp. 397-398:
$$
b(x:y:z) = (yz:xz:xy).
$$
So $b$ is invariant under the natural action of $G_{M/L} := {\rm Gal}(M/L)$.

But the composition $d := (\pi_Y)_M \circ (\pi_X)_M^{-1}$ is also 
$G_{M/L}$-invariant. And $b = \psi \circ d \circ \phi^{-1}$, so 
$$
d = \psi^{-1} \circ b \circ \phi.
$$

Now take $\sigma \in G_{M/L}$. Since ${}^{\sigma} d = d$ and ${}^{\sigma} b = b$,
we have
\begin{align*}
{}^{\sigma}(\psi^{-1} b \phi) &= \psi^{-1} b \phi \\
\psi ({}^{\sigma} \psi^{-1}) &= b \phi ({}^{\sigma} \phi^{-1}) b^{-1} \quad (1)
\end{align*}

Recall that the correspondence between Severi-Brauer varieties with 
points in $M$ and central simple algebras split by $M$ is via the 
cohomology group $H^1(G_{M/L},PGL_3(M))$. The cocycles associated to
$X$ and $Y$ are precisely $\eta_{\sigma} := \phi({}^{\sigma} \phi^{-1})$ and 
$\xi_{\sigma} := \psi({}^{\sigma} \psi^{-1})$. So (1) translates to:
$$
\xi_{\sigma} = b \eta_{\sigma} b^{-1}. \quad (2)
$$

Now, for any $\sigma \in G_{M/L}$, the set 
$\{ (1:0:0), (0:1:0), (0:0:1) \}$ is stable under $\eta_{\sigma}$ and 
$\xi_{\sigma}$, considered as automorphisms of ${\mathbb P}^2$. So $\eta_{\sigma}$
and $\xi_{\sigma}$ land in the subgroup
$$
H := \{ A \in PGL_3(M) : \text{each row and column of $A$ has 
exactly one nonzero entry} \}.
$$

(For simplicity, we often abuse notation and identify elements of 
$PGL_3(M)$ with representative matrices in $GL_3(M)$.)

Clearly $H$ decomposes as a semi-direct product
$$
H = H_D \rtimes H_P,
$$
where $H_D$ is the subgroup of diagonal matrices and $H_P$ the subgroup of 
permutation matrices in $PGL_3(M)$. In particular, any matrix in $H$ can 
be written $A = A_D A_P$, where $A_D$ is diagonal and $A_P$ is a 
permutation matrix.

Conjugation by $b$ sends $A_D$ to $A_D^{-1}$, and sends $A_P$ to $A_P$. Note that
$$
b A_D A_P b^{-1} = A_D^{-1}A_P = \left( {}^t A_D^{-1} \right) \left( {}^t A_P^{-1} 
\right),
$$
or, in other words, conjugation by $b$ is the same as applying the operator
$A \mapsto {}^tA^{-1}$ on $H$.

Next we prove a lemma about this operator.

\begin{lemma} \label{lemma1}
Let $c \colon G_{M/L} \to PGL_3(M)$ be a cocycle corresponding to a central simple 
algebra $x$. Precomposing $c$ with the map $A \mapsto {}^tA^{-1}$ on 
$PGL_3(M)$ sends the class of $x$ to the inverse of the class of $x$ in 
${\rm Br}(L)$.

\end{lemma}

{\em Proof of lemma:} By definition of the correspondence between $c$ and $x$, $c$
is constructed by the following formula: there is some isomorphism 
$\alpha \colon x \otimes_k {\overline k} \to M_3({\overline k})$ such that for any 
$\sigma \in G_{M/L}$ and $A \in M_3({\overline k})$,
$$
c_{\sigma} A c_{\sigma}^{-1} = \alpha ({}^{\sigma} \alpha^{-1})(A).
$$

If we let $\beta(A) = {}^t \alpha(A)$, then $\beta$ is an 
isomorphism $x^{\rm op} \otimes {\overline k} \to M_3({\overline k})$. And
$$
\beta({}^{\sigma} \beta^{-1})(A) = {}^t ((\alpha ({}^{\sigma} \alpha^{-1}))({}^t A))
= {}^t(c_{\sigma} ({}^t A) c_{\sigma}^{-1}) = {}^t c_{\sigma}^{-1} A \, ({}^t 
c_{\sigma}).
$$
So the cocycle corresponding to $x^{\rm op}$ and $\beta$ is precisely the cocycle 
obtained by precomposing $c$ with the map $A \mapsto {}^tA^{-1}$. This proves the 
lemma.
$\Box$

\smallskip

Therefore, by the lemma applied to equation (2), $\eta_{\sigma}$ and $\xi_{\sigma}$ 
correspond to inverse classes in Br$(L)$, i.e. $y$ is Brauer-equivalent to 
$x^{\rm op}$. Since $x$ and $y$
are 9-dimensional, they are either both isomorphic to $M_3(L)$ or both division 
algebras. In the first case, $y \cong x^{\rm op}$ trivially, and in the second 
case $y$ is Brauer-equivalent to $x^{\rm op}$, and two division algebras which are 
Brauer-equivalent are isomorphic. This proves the theorem. $\blacksquare$

\medskip

From now on, we will denote by $X^{\rm op}$ the variety which corresponds to the 
central simple algebra opposite to the one corresponding to $X$; $X^{\rm op}$ is
unique up to isomorphism.

\section{Automorphisms of Severi-Brauer surfaces}

The next result we need is about the action of the automorphism group of a
Severi-Brauer surface $X$ on sets of three non-collinear points.

\begin{theorem} \label{theorem2} Let $X$ be a $2$-dimensional Severi-Brauer variety 
over a number field
$L$, equipped with two $G_L$-stable sets 
$P=\{ P_1,P_2,P_3 \}$ and $Q = \{ Q_1,Q_2,Q_3 \}$ of non-collinear points.
If $\xi: P \to Q$ is an isomorphism of $L$-varieties, $\xi$ can be
extended to an automorphism $\alpha \in {\rm Aut}_L(X)$. 
\end{theorem}

{\em Proof:} Let $M$ be the smallest Galois extension of $L$ over which the points
in $P$ (and $Q$) are all individually defined. Let $G_{M/L} = {\rm Gal}(M/L)$, as 
above. Since $X(M) \ne \emptyset$, we
have an isomorphism $\phi: X_M \to {\mathbb P}^2_M$. As before, since Aut
${\mathbb P}^2_M$ acts transitively on sets of three points in general
position, we may assume that $P_1, P_2, P_3$ go to whatever three
non-collinear points we want. The following easy lemma provides those points:

\begin{lemma} \label{lemma2}
Given a $G_{M/L}$-set $Z$ of order 3, we can find a set $R$ of three non-collinear 
points in ${\mathbb P}^2_M$ such that $R$ and $Z$ are isomorphic as $G_{M/L}$-sets.
\end{lemma}

{\em Proof of lemma:} First note that
we can immediately find a set of three distinct points $\{ \alpha_1, \alpha_2, \alpha_3 
\} \subseteq {\mathbb A}^1(M)$ which is invariant under $G_{M/L}$ and has the desired 
structure as a $G_{M/L}$-set. Let
$R_i = (1 \colon \alpha_i \colon \alpha_i^2)$. Then the $R_i$ are
non-collinear, and the $G_{M/L}$-action on the $R_i$ is the same as the action on the 
$\alpha_i$, which is what we wanted. $\Box$

\smallskip

Applying the lemma with $Z = P$, we obtain a set $R$ of points with the same 
$G_{M/L}$-action as the one on $P$. So set $\phi(P_i) = R_i$ for $i=1,2,3$. We can also
construct an isomorphism $\psi: X_M \to {\mathbb P}^2_M$ such
that $\psi(Q_i) = R_i$ for $i=1,2,3$. Make cocycles $\eta_{\sigma}
= \phi({}^{\sigma}\phi^{-1})$ and $\xi_{\sigma} =
\psi({}^{\sigma}\psi^{-1})$. We know that $\eta_{\sigma}$ and
$\xi_{\sigma}$ are cohomologous in $H^1(G_{M/L},PGL_3(M))$, since they both
correspond to the same Severi-Brauer variety, and this cohomology group 
parameterizes Severi-Brauer varieties split by $M$. But in fact, from the 
construction of $\phi$ and $\psi$ and the fact that $P$, $Q$, and $R$ have the same
$G_{M/L}$-actions, we see that $\eta_{\sigma}$ and $\xi_{\sigma}$ 
can be viewed as cocycles in $Z^1(G_L,T)$, where 
$$ 
T = \{ A \in PGL_3(M) : A(R_i) = R_i \ {\rm for} \ i = 1,2,3 \}. 
$$

We will need to prove the following

\begin{lemma} \label{lemma3} The natural map $i : H^1(G_{M/L},T) \to 
H^1(G_{M/L},PGL_3(M))$ is injective.
\end{lemma}

After the lemma is proved, we will conclude that $\eta_{\sigma}$ and 
$\xi_{\sigma}$ are cohomologous via a coboundary with image in $T$, i.e.
$$
\xi_{\sigma} = B \eta_{\sigma} ({}^{\sigma} B^{-1})
$$
with $B \in T$, so that
\begin{align*}
\psi({}^{\sigma}\psi^{-1}) &= B \phi({}^{\sigma}\phi^{-1})({}^{\sigma}B^{-1}) \\
{}^{\sigma}(\psi^{-1} B \phi) &= \psi^{-1} B \phi \\
\end{align*}
for all $\sigma \in G_L$. So $\psi^{-1} B \phi$ descends to an $L$-automorphism 
which extends $\varphi$.

{\em Proof of lemma:} First, let $U$ be the set of matrices $B \in GL_3(M)$ such 
that the coordinate vectors in $M^3$ representing the $R_i$ are eigenvectors of 
$B$. Then we get the following commutative diagram:

\[ \xymatrix{ 
1 \ar[r] & M^* \ar@{=}[d] \ar[r] & U \ar[d] \ar[r] & T \ar[d] \ar[r] & 1 \\
1 \ar[r] & M^* \ar[r] & GL_3(M) \ar[r] & PGL_3(M) \ar[r] & 1
} \]

Note that $U$ is abelian (indeed, it is clearly conjugate to the subgroup of 
$GL_3(M)$
consisting of the invertible diagonal matrices). So we can pass to the long exact
sequence of cohomology associated to this short exact sequence, part of which is

\[ \xymatrix{
H^1(G_{M/L},U) \ar[r] & H^1(G_{M/L},T) \ar[d]^i \ar[r] & H^2(G_{M/L},M^*) 
\ar@{=}[d] \\
& H^1(G_{M/L},PGL_3(M)) \ar[r] & H^2(G_{M/L},M^*) \\ } \]
 
So if we can show that $H^1(G_{M/L},U) = 0$, we'll have that the map at the top 
of the square is injective, which will imply that $i$ is injective.

First note that $U =CK_DC^{-1}$, where $K_D$ is the subgroup of diagonal matrices 
of $GL_3(M)$, and $C$ is a change-of-basis matrix. So we have a group isomorphism 
$U \to M^* \times M^* \times M^*$ sending $CAC^{-1} \to (A_{11}, A_{22}, A_{33})$.  

Let $R = {\rm Spec} \, E$ as an $L$-variety. $E$ is a 
three-dimensional $L$-algebra, and there are three distinct maps $M 
\otimes_L E \to M$ corresponding to the three elements of $R(M)$. Then 
the group homomorphism $M \otimes_L E \to M \times M \times M$ made from 
these three maps is an isomorphism of rings. Passing to the unit groups 
of both rings gives a group isomorphism $(M \otimes_L E)^* \to M^* 
\times M^* \times M^*$. It is easy to check that the composition $U \to 
M^* \times M^* \times M^* \to (M \otimes_L E)^*$ actually commutes with 
the action of $G_{M/L}$ on both sides. 

Indeed, another way to see this composition is as the realization of $U$ as the
automorphism group of the line bundle over $R_M$ corresponding to the invertible
sheaf ${\mathcal O}_{R_M}(1)$. This automorphism group is isomorphic to
${\mathcal O}_{R_M}(R_M)^* = (M \otimes_L E)^*$.

But $H^1(G_{M/L},(M \otimes_L E)^*) = 0$ by an extension of Hilbert's Theorem 90
(see \cite{serrelocalfields}, X.1, ex. 2). This proves the lemma. 
$\Box$ $\blacksquare$

\section{Reversing the construction}

Now we prove a result about recovering the Del Pezzo surface $D$ from a suitably
chosen Severi-Brauer variety $X$. 

\begin{theorem} \label{theorem3} Let $k$ be a field and let $L/k$ be a quadratic 
extension with ${\rm 
Gal}(L/k)$ generated by $\sigma$. Suppose $X$ is a Severi-Brauer variety over $L$ such 
that $X$ and ${}^{\sigma}X$ correspond to opposite central simple algebras, and 
suppose we are given a $G_L$-stable set of non-collinear points 
$P := \{ P_1, P_2, P_3 \} \subseteq X({\overline k})$. Then:

(i) The variety $S$ we constructed in Proposition \ref{prop3} can be descended to a Del 
Pezzo surface of degree 6 over $k$.

(ii) If we relax the requirements on the above set of data to let $L$ be an
\'etale algebra of degree 2 over $k$, then every Del Pezzo surface of degree 6 over 
$k$ can be constructed in this way.
\end{theorem}

{\em Proof of theorem:} In Proposition \ref{prop3} we can take $Y = {}^{\sigma}X$,
and by Theorem \ref{theorem2} we can assume that the set $Q$ of blown-up points on $Y$ 
is actually ${}^{\sigma}P$. So we have blowing-down maps $S \to X$ and $S
\to {}^{\sigma}X$ as in the proposition, hence a map $\varphi: S \to X 
\times {}^{\sigma}X$. Now we prove 

\begin{lemma} \label{lemma4} $\varphi$ is a closed immersion. \end{lemma}

{\em Proof of lemma:} It is equivalent to show that $\varphi_{\overline k} \colon 
S_{\overline k} \to {\mathbb P}^2_{\overline k} \times {\mathbb 
P}^2_{\overline k}$ is a closed immersion. So $\varphi_{\overline 
k}$ is the map which takes ${\mathbb P}^2$ blown up at three points 
and blows down each skew triple of exceptional curves in turn. This 
description of the map makes it clear that it is injective as a map of 
sets, and since blowups of projective schemes are projective, 
$\varphi_{\overline k}$ is projective; so $\varphi$ is projective and 
thus $\varphi$ is a homeomorphism onto its image, a closed subset of 
${\mathbb P}^2 \times {\mathbb P}^2$. We now need to check that the map on 
structure sheaves is surjective, which can be checked on the stalks.

What we need to check is that the map ${\mathcal O}_{{\mathbb P}^2 \times {\mathbb
P}^2, (\pi_X(P),\pi_Y(P))} \to {\mathcal O}_{S,P}$ induced by $\varphi$ is 
surjective for all $P \in S({\overline k})$. (For convenience, we assume for the
remainder of the lemma that everything is over $\overline k$ and drop subscripts.)
If $P$ lies on at most one of the 
exceptional lines, then one of the projections $p_S: S \to {\mathbb P}^2 \times 
{\mathbb P}^2 \to {\mathbb P}^2$ restricts to an isomorphism of an open subset of 
$S$ containing $P$ (namely, $S$ minus a skew triple of exceptional lines not 
containing $P$) onto its image. Thus 
$$
{\mathcal O}_{{\mathbb P}^2,p_S(P)} \to {\mathcal O}_{{\mathbb P}^2 \times {\mathbb 
P}^2,\varphi(P)} \to {\mathcal O}_{S,P}
$$
is surjective, and so the latter map must be as well.

Now suppose $P$ is one of the six points which lies on two exceptional lines. As in 
\cite{hartshorne1977}, p. 152, it is enough to check that the map $m_{{\mathbb
P}^2 \times {\mathbb P}^2,\varphi(P)} \to m_{S,P}/m_{S,P}^2$ is surjective. Around $P$,
$S$ just looks like the blowup of ${\mathbb A}^2$ at a point, and so 
$m_{S,P}/m_{S,P}^2$ is two-dimensional, with generators which cut out the two 
exceptional lines going through $P$. Each of these two generators comes from exactly
one of the maps $m_{{\mathbb P}^2,p_S(P)} \to m_{S,P}/m_{S,P}^2$ (whichever one 
does not collapse the line that that generator cuts out). So this implies the
surjectivity of the map we want. $\Box$

\smallskip

Now since $\varphi$ is a closed immersion, it gives
an isomorphism of $S$ onto its image, which must be the graph of the
birational map $b_1: X \dashrightarrow S \to {}^{\sigma}X$ that blows up
the $P_i$ and then blows down the ``other" three lines into the
$Q_i$. (The graph of a birational map $b$ is the {\em closure} of the set of points
$\{(a,b(a)) | a \in \text{domain$(b)$} \}$.)

The same analysis shows that ${}^{\sigma}S$ is isomorphic to the graph of the 
birational map $b_2: {}^{\sigma}X \dashrightarrow S \to X$ blowing up the 
$Q_i$ and then blowing down the ``other" three lines into the $P_i$. Since 
$b_1$ and $b_2$ are inverses by construction, the identification of 
$X \times {}^{\sigma}X$ with ${}^{\sigma}X \times X$ by changing the order of the 
factors induces a map $f_{\sigma}$ from the graph of $b_1$ to the graph of $b_2$. Then
we obtain the following commutative diagram:

\[ \xymatrix{
S \ar[r] \ar[d]^{f_{\sigma}} & X \times {}^{\sigma}X 
\ar[d] \\
{}^{\sigma}S \ar[r] \ar[d]^{{}^{\sigma}f_{\sigma}} & {}^{\sigma}X 
\times X \ar[d] \\
S \ar[r] & X \times {}^{\sigma}X \\
} \]

where the maps on the right are the isomorphisms arising from switching the factors. 
This shows that ${}^{\sigma}f_{\sigma} \circ f_{\sigma}$ is the identity, since the 
composition of maps on the right side is the identity. Therefore $f_{\sigma}$ gives 
descent data for $S$, so that it can be descended (again, as in \cite{weil1956}) to a 
$k$-variety. This proves statement (i).

As for statement (ii), this merely summarizes what we already know. Given a Del Pezzo
surface $D$, if the field given in Proposition \ref{prop2} was quadratic, we can take
it to be $L$, and if that field was $k$, we can take $L = k \times k$.
Then the Severi-Brauer varieties $X$ and $Y$ associated with $D$ in
Proposition \ref{prop2} satisfy $X^{\rm op} = Y = {}^{\sigma}X$ in either case. 
$\blacksquare$

\section{Other results and conclusions}

One remark that should be made about Theorem \ref{theorem3} is that the condition 
${}^{\sigma}X = X^{\rm op}$ is not a very restrictive one. When $L = k \times k$, 
the varieties $X$ satisfying the condition are generated by starting with any 
Severi-Brauer surface $X'/k$ and then letting $X$ be the disjoint union of $X'$ and 
$(X')^{\rm op}$. 

When $L$ is a quadratic field extension of $k$, we can view the Galois group 
generator $\sigma$ as a linear automorphism of Br $L$. Suppose $x$ is a class of
order $3$ in Br $L$. Then $x$ can be written as $2(x + {}^{\sigma}x) + 
2(x - {}^{\sigma} x)$, so 
$$
({\rm Br} \ L)[3] = ({\rm Br} \ L)^{G_{L/k}}[3] \oplus W,
$$
where $W$ is the set of classes of algebras corresponding to varieties $X$ 
satisfying the condition ${}^{\sigma} X = X^{\rm op}$.

In fact, a little more can be said: the spectral sequence 
$$
E_2^{p,q} := H^p(G_{L/k},H^q(G_L,{\overline L}^*)) \Rightarrow H^{p+q}(G_k,
{\overline k}^*)
$$
yields the usual exact sequence
$$
0 \to E_2^{2,0} \to {\rm Br} \ k \to E_2^{0,2} \to E_2^{3,0},
$$
but $E_2^{0,2} = ({\rm Br} \ L)^{G_{L/k}}$ and $E_2^{3,0} = E_2^{1,0} = 0$ by
Hilbert's Theorem 90 and the fact that $L/k$ is cyclic. Since multiplication-by-3
is the identity on the $2$-torsion group $E_2^{2,0}$, the natural map
$$
({\rm Br} \ k)[3] \to ({\rm Br} \ L)^{G_{L/k}}[3]
$$
is an isomorphism, so that $({\rm Br} \ L)[3] \cong ({\rm Br} \ k)[3] \oplus W$.

Finally, we simultaneously sum up the results we have established and include the
proof of one last remark:

\begin{theorem} \label{theorem4} Giving a Del Pezzo surface of degree $6$ over a field 
$k$ of characteristic zero is equivalent to giving the following data:

\begin{enumerate}

\item an \'etale algebra $L$ of degree $2$ over $k$

\item a Severi-Brauer variety $X$ of dimension $2$ over $L$ such that 
${}^{\sigma} X = X^{\rm op}$, where $\sigma$ generates Gal$(L/k)$

\item a subscheme $P$ of $X$ consisting of three geometric non-collinear points

\end{enumerate}

Moreover, two Del Pezzo surfaces $S_i$ corresponding to $L_i$, $X_i$, and $P_i$ 
$(i=1,2)$ are isomorphic if and only if $L_1 \cong L_2$ and there is an isomorphism 
$X_1 \to X_2$ such that $P_1$ maps isomorphically onto $P_2$. 

\end{theorem}

{\em Proof of theorem:} All we need to check is the last statement. For the ``if"
direction, this simply follows from the description of $S_i$ given in the
proof of Theorem \ref{theorem3}, as the graph of a birational map constructed in 
terms of $X_i$, ${}^{\sigma} X_i$, and $P_i$. For the ``only if" direction, note that
we constructed
the objects $L_i$, $X_i$, and $P_i$ intrinsically from $S_i$. If $S_1$ and $S_2$ are
isomorphic, we naturally get the isomorphisms given in the statement of the theorem.
(The only choice we made was between $X$ and ${}^{\sigma}X$, but these varieties are
isomorphic over $L$, and, as noted before, Theorem \ref{theorem2} implies that we can 
make the isomorphism send $P$ to ${}^{\sigma} P$.) $\blacksquare$

\section*{Acknowledgements}

I thank my advisor Bjorn Poonen, who posed and helped answer many of the questions
raised by this paper.

\bibliographystyle{alpha}
\bibliography{pat}

\end{document}